\magnification=1200

%%%%%%%%%%%%%%%%%%%%%%%%%%%%%%%%%%%%%%%%%%%%%%%%%%%%
% typoref.tex. V : January 18, 2000. 
% Author : Anthony PHAN
% Warning : syntaxe +- LaTeX 
% Sources :
% T. Lachand--Robert, ``La Ma\^\i trise de \TeX'',
% R\'ef\'erences crois\'ees;
% latex.ltx's sources;
% and of course the \TeX book.
%%%%%%%%%%%%%%%%%%%%%%%%%%%%%%%%%%%%%%%%%%%%%%%%%%%%%
%
\catcode`@=11
%
% style (look at the behavior of \item dans \bibitem too,
% and at one ,\  in \re@dreferenceslist)
% Feel free to change: 	\bibn@me (title like ``R\'ef\'erences'')
%			\bibliographym@rk (general style)
%
\def\bibn@me{R\'ef\'erences}
\def\bibliographym@rk{\centerline{{\sc\bibn@me}}
	\sectionmark\section{\ignorespaces}{\unskip\bibn@me}
	\bigbreak\bgroup
	\ifx\ninepoint\undefined\relax\else\ninepoint\fi}
%
% Beware of the \bgroup: it will be closed by \endthebibliography
%
% \refsp@ce is the spacing command that appens between multiple
% references.
%
\let\refsp@ce=\ 
\let\bibleftm@rk=[
\let\bibrightm@rk=]
%
% if you want more space between brackets...
%\let\refsp@ce=\thinspace
%\def\bibleftm@rk{[\thinspace}
%\def\bibrightm@rk{\thinspace]}
%
% frenchy stuff
%
\def\numero{n\raise.82ex\hbox{$\fam0\scriptscriptstyle o$}~\ignorespaces}
%
% new variables
%
\newcount\equationc@unt
\newcount\bibc@unt
\newif\ifref@changes\ref@changesfalse
\newif\ifpageref@changes\ref@changesfalse
\newif\ifbib@changes\bib@changesfalse
\newif\ifref@undefined\ref@undefinedfalse
\newif\ifpageref@undefined\ref@undefinedfalse
\newif\ifbib@undefined\bib@undefinedfalse
\newwrite\@auxout
%
% mark an equation
%
\def\eqnum{\global\advance\equationc@unt by 1%
\edef\lastref{\number\equationc@unt}%
\eqno{(\lastref)}}
%
% One can reference anything, just copy the former macro
% and use it so: \machin \label{truc}
% In machin you would have defined \lastref by some number
% or any text.
%
% References macros
%
% The next macros are the core of \ref and \cite commands.
% Its first argument may be ref, pageref or bib.
%
% It is too tricky to be explained.
% (It is a bit recursive.)
% It allows using \cite or \ref or ...
% with arbitrary many arguments,
% for instance:
% \cite{knuth1,knuth2,ma pomme}
%
% First argument is always ref, pageref or bib.
%
\def\re@dreferences#1#2{{%
	\re@dreferenceslist{#1}#2,\undefined\@@}}
\def\re@dreferenceslist#1#2,#3\@@{\def\next{#2}%
	\expandafter\ifx\csname#1@@\meaning\next\endcsname\relax
	??\immediate\write16
	{Warning, #1-reference "\next" on page \the\pageno\space
	is undefined.}%
	\global\csname#1@undefinedtrue\endcsname
	\else\csname#1@@\meaning\next\endcsname\fi
	\ifx#3\undefined\relax
	\else,\refsp@ce\re@dreferenceslist{#1}#3\@@\fi}
%
% notice that the former ``,\refsp@ce'' will separate
% multiple arguments. But beware of spaces
% while defining a reference or calling for it!
%
% tricky thing: \newlabel has two arguments
% {labelname}{{\lastref}{\pageref}}
% The second argument is read as two arguments
% by \newl@bel. This was necessary to get
% a jobname.aux containing the same syntax
% LaTeX would produce and use.
%
\def\newlabel#1#2{{\def\next{#1}\newl@bel#2}}
\def\newl@bel#1#2{%
	\expandafter\xdef\csname ref@@\meaning\next\endcsname{#1}%
	\expandafter\xdef\csname pageref@@\meaning\next\endcsname{#2}}
\def\label#1{{%
	\toks0={#1}\message{ref(\lastref) \the\toks0,}%
	\ignorespaces\immediate\write\@auxout%
	{\noexpand\newlabel{\the\toks0}{{\lastref}{\the\pageno}}}%
	\def\next{#1}%
	\expandafter\ifx\csname ref@@\meaning\next\endcsname\lastref%
	\else\global\ref@changestrue\fi%
	\newlabel{#1}{{\lastref}{\the\pageno}}}}
\def\ref#1{\re@dreferences{ref}{#1}}
\def\pageref#1{\re@dreferences{pageref}{#1}}
%
% bibliography macros
%
\def\bibcite#1#2{{\def\next{#1}%
	\expandafter\xdef\csname bib@@\meaning\next\endcsname{#2}}}
\def\cite#1{\bibleftm@rk\re@dreferences{bib}{#1}\bibrightm@rk}
%
% The argument of \beginthebibliography
% is any sequence of numerals which will represent
% the maximum \item's length. If you have less than 9
% \bibitem's, this argument may be {any numeral}.
% if you have between 100 and 999 \bibitem's
% this argument may be {any three numerals},
% and so on.
%
\def\beginthebibliography#1{\bibliographym@rk
	\setbox0\hbox{\bibleftm@rk#1\bibrightm@rk\enspace}
	\parindent=\wd0
	\global\bibc@unt=0
	\def\bibitem##1{\global\advance\bibc@unt by 1
		\edef\lastref{\number\bibc@unt}
		{\toks0={##1}
		\message{bib[\lastref] \the\toks0,}%
		\immediate\write\@auxout
		{\noexpand\bibcite{\the\toks0}{\lastref}}}
		\def\next{##1}%
		\expandafter\ifx
		\csname bib@@\meaning\next\endcsname\lastref
		\else\global\bib@changestrue\fi%
		\bibcite{##1}{\lastref}
		\medbreak
		\item{\hfill\bibleftm@rk\lastref\bibrightm@rk}%
		}
	}
\def\endthebibliography{\egroup\par}
%
% THE NEXT MACRO MUST BE INCLUDED
% IN THE \BYE COMMAND. FOR INSTANCE:
%
% \catcode`@=11
% \outer\def\bye{\@closeaux
% 	\par\vfill\supereject\end}
% \catcode`@=12
%
\def\@closeaux{\closeout\@auxout
	\ifref@changes\immediate\write16
	{Warning, changes in references.}\fi
	\ifpageref@changes\immediate\write16
	{Warning, changes in page references.}\fi
	\ifbib@changes\immediate\write16
	{Warning, changes in bibliography.}\fi
	\ifref@undefined\immediate\write16
	{Warning, references undefined.}\fi
	\ifpageref@undefined\immediate\write16
	{Warning, page references undefined.}\fi
	\ifbib@undefined\immediate\write16
	{Warning, citations undefined.}\fi}
%
% initialization of jobname.aux
%
\immediate\openin\@auxout=\jobname.aux
\ifeof\@auxout \immediate\write16
  {Creating file \jobname.aux}
\immediate\closein\@auxout
\immediate\openout\@auxout=\jobname.aux
\immediate\write\@auxout {\relax}%
\immediate\closeout\@auxout
\else\immediate\closein\@auxout\fi
%
% Let's read this file and open it out
%
\input\jobname.aux
\immediate\openout\@auxout=\jobname.aux
% this file will be closed by \bye.
%
% That's all, folks!
%
\catcode`@=12
%\endinput

%
\catcode`@=11
\def\bibliographym@rk{\bgroup}
%
% \bye est modifie pour la biblio et la table des matieres
%
\outer\def\bye{ 	\par\vfill\supereject\end}

\def\house#1{\setbox1=\hbox{$\,#1\,$}%
\dimen1=\ht1 \advance\dimen1 by 2pt \dimen2=\dp1 \advance\dimen2 by 2pt
\setbox1=\hbox{\vrule height\dimen1 depth\dimen2\box1\vrule}%
\setbox1=\vbox{\hrule\box1}%
\advance\dimen1 by .4pt \ht1=\dimen1
\advance\dimen2 by .4pt \dp1=\dimen2 \box1\relax}

  \def\eps{{\varepsilon}}

  \def\noi{\noindent}

\def\build#1_#2^#3{\mathrel{\mathop{\kern 0pt#1}\limits_{#2}^{#3}}}

\def\date {le\ {\the\day}\ \ifcase\month\or janvier
\or fevrier\or mars\or avril\or mai\or juin\or juillet\or
ao\^ut\or septembre\or octobre\or novembre
\or d\'ecembre\fi\ {\oldstyle\the\year}}

\font\fivegoth=eufm5 \font\sevengoth=eufm7 \font\tengoth=eufm10

\newfam\gothfam \scriptscriptfont\gothfam=\fivegoth
\textfont\gothfam=\tengoth \scriptfont\gothfam=\sevengoth

\def\pro{\noindent {\it Proof. }}

\def\smallsquare{\vbox{\hrule\hbox{\vrule height 1 ex\kern 1 ex\vrule}\hrule}}
\def\cqfd{\hfill \smallsquare\vskip 3mm}

\def\og{\leavevmode\raise.3ex\hbox{$\scriptscriptstyle 
\langle\!\langle\,$}}
\def \fg {\leavevmode\raise.3ex\hbox{$\scriptscriptstyle 
\!\rangle\!\rangle\,\,$}}

\def\rme{{\rm e}}

\def\bfv{{\bf v}}

%%%%%%%%%%%%%%%%%%%%%%%%%%%%%%%%%%%%%%%%%%%%
\centerline{}

\vskip 8mm

\centerline{\bf On the digital representation of smooth numbers}

\vskip 13mm

\centerline{Yann B{\sevenrm UGEAUD} and Hajime K{\sevenrm{ANEKO}} \footnote{}{\rm
}}

{\narrower\narrower
\vskip 15mm

\proclaim Abstract. {
Let $b \ge 2$ be an integer. 
Among other results, we establish, in a quantitative form,  
that any sufficiently large integer which is not a multiple of $b$
cannot have simultaneously only few distinct prime factors and only few nonzero digits in 
its representation in base $b$. 
}

}

\vskip 10mm

\centerline{\bf 1. Introduction and results}

\vskip 5mm

Let $a, b$ be positive, multiplicatively independent integers.
Stewart \cite{Ste80} established that, for every sufficiently large
integer $n$, the representation of $a^n$ in base $b$ 
has more than $(\log n) / (2 \log \log n)$ nonzero digits. His proof rests 
on a subtle application of Baker's theory of linear forms in complex logarithms of 
algebraic numbers. This result addresses a very special case of the 
following general (and left intentionally vague) question, which was introduced and 
discussed in \cite{Bu17}:

\smallskip
{\it Do there exist arbitrarily large integers 
which have only small prime factors and, at the same time, few nonzero digits in their
representation in some integer base?}

\smallskip

The expected answer is {\it no} and a very modest step in this direction 
has been made in \cite{Bu17}, 
by using a combination of estimates for linear forms in complex and $p$-adic logarithms. 
In the present work, we considerably extend Corollary 1.3 of \cite{Bu17} and,
more generally, we show in a quantitative form that the
maximum of the greatest prime factor of an integer $n$ and the number 
of nonzero digits in its representation in a given integer base 
tends to infinity  
as $n$ tends to infinity.

Throughout this note, $b$ always denotes an integer at least equal to $2$. 
Following \cite{Bu17}, for an integer $k \ge 2$, we denote by $(u_j^{(k)})_{j \ge 1}$ 
the sequence, arranged in increasing order, of all positive integers which are 
not divisible by $b$ and have at most $k$ nonzero digits in their representation in base $b$. 
Said differently, $(u_j^{(k)})_{j \ge 1}$ is the ordered sequence composed of the integers 
$1, 2, \ldots , b-1$ and those of the form
$$
d_k b^{n_k} + \cdots + d_2 b^{n_2} + d_1, 
\quad n_k > \cdots > n_2 > 0, \quad  
d_1, \ldots , d_k \in \{0, 1, \ldots , b-1\}, \quad
d_1 d_k \not= 0.
$$
We stress that, for the questions investigated in the present note, 
it is natural to restrict 
our attention to integers not divisible by $b$. 
Obviously, the sequence $(u_j^{(k)})_{j \ge 1}$ depends on $b$, but, for 
shortening the notation, we have decided not to mention this dependence.

Theorem 1.1 of \cite{Bu17} implies that 
the greatest prime factor of $u_j^{(k)}$ tends to infinity as $j$ tends to infinity.
Its proof rests on the Schmidt Subspace Theorem and does not allow us
to derive an estimate for the speed of convergence. 
Such an estimate has been  
established in \cite{Bu17}, 
but only for $k \le 3$. 
Following the proof of our main result (Theorem 1.2 below), we are able
to extend this estimate to arbitrary integers $k$.

For a positive integer $n$, let denote by $P[n]$ its greatest prime factor 
and by $\omega (n)$ the number 
of its distinct prime factors, with the convention that $P[1] = 1$. 
A positive real number $B$ being given, 
a positive integer $n$ is called $B$-smooth if $P[n] \le B$. 
%Throughout this note, denote by $p_j$ the $j$-th prime number for $j \ge 1$.   %%y 

\proclaim Theorem 1.1. 
Let $b \ge 2, k\ge 3$ be integers. 
Let $\eps$ be a positive real number.
Then, there exists an effectively
computable positive number $j_0$, depending 
only on $b,k$ and $\eps$, such that 
$$
P[u_j^{(k)}] > \Bigl({1\over k-2} - \eps\Bigr) \log \log u_j^{(k)} 
\, {\log \log \log u_j^{(k)} \over \log \log \log \log u_j^{(k)}}, \quad
\hbox{for $j > j_0$}.    
$$
In particular, there exists an effectively computable positive integer $n_0$, depending 
only on $b,k$ and $\eps$, such that any integer $n > n_0$ which is 
not divisible by $b$ and is 
$$
\Bigl({1\over k-2} - \eps\Bigr) (\log \log n) {\log \log \log n \over \log \log \log \log n}\hbox{-smooth}
$$
has at least $k+1$ nonzero digits in its $b$-ary representation. 

Taking $k=3$ in Theorem 1.2, we get the second assertion 
of Theorem 1.3 of \cite{Bu17}. 
The main ingredients for the proofs of both theorems are estimates 
for linear forms in complex and $p$-adic logarithms of algebraic numbers. 
The novelty in the present note is a repeated use of estimates for linear forms in $p$-adic 
logarithms, where $p$ is a prime divisor of the base $b$. 
With our new approach, the number $k$ of nonzero digits need not 
to be fixed and can be allowed to depend on $n$, provided that it is rather small 
compared to $n$.

Our main result asserts that, given an integer $b \ge 2$, if 
the integer $n$ is sufficiently large, then its greatest prime factor and the 
number of nonzero digits in its representation in base $b$ cannot be simultaneously small.

\proclaim Theorem 1.2. 
Let $b \ge 2$ and $k \ge 2$ be integers. 
There exist an effectively computable real number $c$, depending at most on $b$, and 
an effectively computable, absolute real number $C$ such that
every sufficiently large positive integer $n$, which is not divisible by $b$ 
and whose representation in base $b$ has $k$ nonzero digits, satisfies
$$
{\log \log n\over k} \le c + \log k + \omega (n) (C + \log \log P[n])
+ \log \log (k \log P[n]).    
$$

Several easy consequences of (the proof of) Theorem 1.2 are pointed out below. 
We extend the definition of the sequences $(u_j^{(k)})_{j \ge 1}$ as follows. 
For a positive real valued function $f$ defined 
over the set of positive integers, we
let $(u_j^{(f)})_{j \ge 1}$ be the sequence, arranged in increasing order, 
of all positive integers $n$ which are not 
divisible by $b$ and have at most $f(n)$ nonzero digits in their representation in base $b$.

\proclaim Theorem 1.3.    
Let $b\geq 2$ be an integer. Let $f$ be a positive real valued function defined 
over the set of positive integers such that
$$
\lim_{u\to + \infty} f(u)= +\infty. 
$$
Assume that there exists a real number $\delta$ satisfying $0<\delta<1$ and 
$$
f(u)\leq (1-\delta) \, {\log\log u\over \log\log\log u},
\eqno (1.1) 
$$
for any sufficiently large $u$, and set 
$$
\Psi_f(u):={\log\log u \over f(u)}, \quad \hbox{for $u \ge 3$}.
$$
Then, for an arbitrary positive real number $\eps$, we have 
$$
P[u_j^{(f)}] > (\delta_0-\eps)\Psi_f \bigl(u_j^{(f)}\bigr) \,
{\log \Psi_f \bigl(u_j^{(f)}\bigr) \over \log \log \Psi_f \bigl(u_j^{(f)}\bigr)},   \eqno (1.2)
$$
for any sufficiently large integer $j$, where 
$$
\delta_0=\sup \left\{ \delta>0 \ :  \ f(u)\leq (1-\delta)\, {\log \log u \over \log \log \log u}
\hbox{ for every large integer } u \right\}. 
$$

We gather in the next statement three immediate consequences
of Theorem 1.3 applied with an appropriate function $f$.  

\proclaim Corollary 1.4. 
Let $b \ge 2$ be an integer. 
There exists an effectively computable positive integer $n_0$, depending 
only on $b$, such that any integer $n > n_0$ which is 
not divisible by $b$ satisfies the following three assertions. 
If $n$ is
$$
{\log \log n \over 2 \log \log \log \log n} \hbox{-smooth},
%$$
\,  \hbox{then $n$ has at least} \, \, 
%$$
\log \log \log n 
$$
nonzero digits in its representation in base $b$. 
If $n$ is
$$
\sqrt{ \log \log n \, {\log \log \log n \over \log \log \log \log n} \,} \hbox{-smooth}, 
%$$
\,  \hbox{then $n$ has at least} \, \, 
%$$
{1 \over 3} \, \sqrt{ \log \log n \, {\log \log \log n \over \log \log \log \log n} \, }
$$ 
nonzero digits in its representation in base $b$. 
If $n$ is
$$
{1 \over 2} \, \log \log \log n \, {\log \log \log \log n \over \log \log \log \log \log n} \hbox{-smooth},
%$$
\,  \hbox{then $n$ has at least} \, \, 
%$$ 
{\log \log n \over 2 \log \log \log n} 
$$
nonzero digits in its representation in base $b$.

Let $S$ be a finite, non-empty set of prime numbers. 
A rational integer is an integral $S$-unit if all its prime factors belong to $S$. 
We deduce from Theorem 1.2 
lower bounds for the number of nonzero digits in the representation
of integral $S$-units in an integer base.

\proclaim Corollary 1.5.  
Let $b\geq 2$ be an integer. Let $S=\{q_1<\cdots<q_s\}$ be a set of $s$ distinct prime numbers. 
Then, for any positive real number $\eps$, 
there exists an effectively computable positive integer $n_0$, 
depending only on $b,S$, and $\eps$, such that any integral $S$-unit $n\geq n_0$ 
which is not divisible by $b$ has more than 
$$
(1-\eps){\log\log n \over\log\log\log n} 
$$
nonzero digits in its representation in base $b$. 

Let $a \ge 2, b \ge 2$ be coprime integers.  %%y
By taking for $S$ the set of prime divisors of $a$, Corollary 1.5 implies 
Stewart's result 
mentioned in the introduction (for the case where $a$ and $b$ are multiplicatively independent 
and not coprime, the proof of Corollary 1.5 can be easily adapted) and both proofs are different. 
Observe, however, that 
Stewart obtained in  \cite{Ste80} a more general result, namely that, for any 
multiplicatively independent positive integers $b$ and $b'$ and any sufficiently 
large integer $n$, the number of nonzero digits in the representation of $n$ in 
base $b$ plus the number of nonzero digits in the representation of $n$ in 
base $b'$ exceeds $(\log\log n) / (2 \log\log \log n)$.

Our results are established in Section 3, by means of 
lower estimates for linear forms in logarithms gathered in Section 2. 
We postpone to Section 4 comments and remarks.

\vskip 5mm

\centerline{\bf 2. Lower estimates for linear forms in logarithms}

\vskip 5mm

The first assertion of Theorem 2.1 is an immediate 
consequence of a theorem of Matveev \cite{Matv00}. 
The second one is a slight simplification 
of the estimate given on page 190 of Yu's paper \cite{Yu07}. 
For a prime number $p$ and a nonzero rational number $z$ we denote by $v_p(z)$ the exponent 
of $p$ in the decomposition of $z$ in product of prime factors.

\proclaim Theorem 2.1.  
Let $n \ge 2$ be an integer. 
Let $x_1/y_1, \ldots, x_n/y_n$ be nonzero rational numbers. 
Let $b_1, \ldots, b_n$ be integers such that $(x_1/y_1)^{b_1} \cdots (x_n/y_n)^{b_n} \not= 1$. 
Let $A_1, \ldots, A_n$ be real numbers with
$$
A_i \ge \max\{|x_i|, |y_i|, \rme\}, \quad 1\le i \le n.
$$
Set
$
B = \max\{ 3, |b_1|, \ldots , |b_n| \}.
$
Then, we have
$$
\log \Bigl| \Bigl( {x_1 \over y_1} \Bigr)^{b_1} \cdots \Bigl( {x_n \over y_n} \Bigr)^{b_n}  - 1 \Bigr|  
> - 8 \times 30^{n+3} \, n^{9/2}  \,  \log (\rme B) \, 
\log A_1 \cdots \log A_n.      \eqno (2.1)  
$$
Let $p$ be a prime number. Then, we have
$$
\eqalign{
v_p \Bigl( \Bigl( {x_1 \over y_1} \Bigr)^{b_1} \cdots \Bigl( {x_n \over y_n} \Bigr)^{b_n}  
& - 1 \Bigr) < \cr
& (16 \rme )^{2(n+1)} n^{5/2} (\log (2n))^2
\,  {p \over (\log p)^2} \log A_1 \cdots \log A_n \log B.  \cr}   \eqno (2.2)  
$$

\vskip 5mm

\centerline{\bf 3. Proofs}

\vskip 5mm

Below, the constants $c_1, c_2, \ldots$ are effectively computable and depend at most 
on $b$ and the constants $C_1, C_2, \ldots$ are absolute and effectively computable. 
Let $N$ be a positive integer and $k$ the number of nonzero digits in its representation 
in base $b$. We assume that $b$ does not divide $N$, thus $k \ge 2$ and we write 
$$
N =: d_k b^{n_k} + \cdots + d_2 b^{n_2} + d_1 b^{n_1}, 
$$
where
$$
n_k > \cdots > n_2 > n_1=0, \quad
d_1, \ldots , d_k \in \{1, \ldots , b-1\}.
$$
Let $q_1, \ldots , q_s$ denote distinct prime numbers written in increasing order such that
there exist non-negative integers $r_1,\ldots,r_s$ with 
$$
N=q_1^{r_1}\cdots q_s^{r_s}.
$$
Observe that 
$$
b^{n_k} \le N < b^{n_k + 1}.  \eqno (3.1)
$$

\proclaim Lemma 3.1. 
Keep the above notation and set
$ 
k^* := \max\{k-2, 1\}.
$ 
If 
$$
\log N \ge 2 (\log b) \, \Bigl( {8 \log b \over \log 2} \Bigr)^k,  \eqno (3.2)
$$
then we have
$$
n_k \leq \Bigl(c_1 C_1^{s} k^*
\Bigl( \prod_{i=1}^s \log q_i \Bigr) \, \log (k \log q_s)\Bigr)^{k^*} .
\eqno (3.3)
$$

\pro 
First we assume that $n_k\geq 2 n_{k-1}$. This covers the case $k = 2$. 
Since 
$$
\eqalign{
\Lambda_a:=
\Bigl| \Bigl(\prod_{i=1}^s q_i^{r_i}\Bigr)d_k^{-1} b^{-n_k} -1\Bigr| & =
d_k^{-1} b^{-n_k}\sum_{h=1}^{k-1} d_h b^{n_h} \cr
&\leq b^{1+n_{k-1}-n_k}\leq b^{-(n_k-2)/2}, \cr
}
$$
we get 
$$
\log \Lambda_a  \le -  \Bigl({n_k \over 2} - 1 \Bigr) \, \log b. \eqno (3.4) 
$$
Sicne $r_j \log q_j\le (n_k+1)\log b$ for $j=1,\ldots,s$, we deduce from (2.1) that 
$$
\log \Lambda_a\geq -c_2 C_2^s (\log q_1)\cdots (\log q_s) (\log n_k).
\eqno (3.5)
$$
Combining (3.4) and (3.5), we obtain 
$$
n_k\leq c_3 C_3^s \Bigl( \prod_{i=1}^s \log q_i \Bigr) \, (\log \log q_s),
$$
which implies (3.3). 

Now, we assume that $n_k<2n_{k-1}$. 
In particular, we have $k \ge 3$. 
If there exists an integer $j$ with $1\leq j\leq k-3$ and 
$n_{1+j}\geq n_k^{j/(k-2)}$, then put 
$$
\ell :=\min \{j : 1\leq j\leq k-3, \ n_{1+j}\geq n_k^{j/(k-2)}\}.
$$
Otherwise, set $\ell :=k-2$. We see that 
$$
n_{\ell+1}\geq {1 \over 2} n_k^{\ell /(k-2)}
\quad \hbox{and} \quad 
n_\ell \leq n_k^{(\ell -1)/(k-2)}.
\eqno (3.6)
$$
Let $p$ be the smallest prime divisor of $b$. 
Put
$$
\Lambda_u:=
\Bigl(\prod_{i=1}^s q_i^{r_i}\Bigr)\Bigl(\sum_{h=1}^\ell d_h b^{n_h}\Bigr)^{-1}-1
=\Bigl(\sum_{h=\ell +1}^k d_h b^{n_h}\Bigr)\Bigl(\sum_{h=1}^\ell d_h b^{n_h}\Bigr)^{-1}.
$$
%Suppose that $j$ is sufficiently large in terms of $b,k$ satisfying $n_k^{1/(k-2)}\geq (8\log b)/(\log 2)$. 
We get by (3.6), (3.1), and (3.2) that 
$$
\eqalign{
v_p(\Lambda_u)&\geq n_{\ell+1}-{\log b^{1+n_\ell} \over \log p}\cr
&\geq {1 \over 2} n_k^{\ell/(k-2)}-
\left(
1+n_k^{(\ell-1)/(k-2)}
\right){\log b \over \log 2} \cr
& \ge {1 \over 2} n_k^{\ell/(k-2)} - 2 n_k^{\ell/(k-2)} \, {\log b \over n_k^{1/(k-2)} \log 2}
\geq {1\over 4}n_k^{\ell/(k-2)}. 
\cr} \eqno (3.7)
$$
We deduce from (2.2) and (3.6) that  
$$
v_p(\Lambda_u)\le c_4 C_4^s (\log q_1)\cdots (\log q_s) n_k^{(\ell-1)/(k-2)} \log n_k. 
\eqno (3.8)
$$
By combining (3.7) and (3.8), we get
$$
n_k^{1/(k-2)}\leq c_5 C_5^s (\log q_1)\cdots (\log q_s) (k-2)\log \bigl(n_k^{1/(k-2)}\bigr),
$$
which implies (3.3) and completes the proof of Lemma 3.1.
\cqfd 

\medskip

\noindent {\it Proof of Theorem 1.1.}

\medskip

We keep the above notation. In particular, $N$ denotes an integer not divisible by $b$ 
and with exactly $k$ nonzero digits in its representation in base $b$. 
In view of \cite{Bu17}, we assume that $k \ge 3$, thus $k^* = k-2$. 
Note that (3.2) holds if $N$ is large enough. Then,
we deduce from (3.1) and (3.3) that 
$$
{\log \log N \over k-2} \le c_6 +  C_6 s +  \log(k-2)
+  \sum_{i=1}^s \log \log q_i  + \log \log (k \log q_s).       \eqno (3.9) 
$$
In particular, denoting by $p_j$ the $j$-th prime number for $j \ge 1$ 
and defining $s$ by $p_s = P[N]$, inequality (3.9) applied 
with $q_j = p_j$ for $j = 1, \ldots , s$ shows that 
$$
{\log \log N \over k-2}  \le  c_6 +  C_6 s
+  \log(k-2)+ s \log \log P[N]  + \log \log (k \log P[N]).    \eqno (3.10)
$$
%For a positive integer $n$, let $\omega (n)$ denote the number 
%of its distinct prime factors. By the Prime Number Theorem, for every positive real
%number $\eps$, there exists an effectively computable integer $n_0 (\eps)$ such that 
%Define $s$ by $p_s = P[n]$ for an integer $n\geq 2$. 
Let $\eps$ be a positive real number. 
By the Prime Number Theorem,  
there exists an effectively computable integer $s_0 (\eps)$, depending only on $\eps$, such that, 
if $s\geq s_0 (\eps)$, then 
$$
s<  (1 + \eps) \, {P [N] \over \log P[N]}.    \eqno (3.11)
$$
For $s < s_0 (\eps/(2k-4))$, we derive from (3.10) an upper bound for %%y 
$N$ in terms of $b,k$ and $\eps$.
For $s\geq s_0 (\eps/(2k-4))$, it follows from 
(3.11) and the Prime Number Theorem that
$$
\log \log N \le (k-2+ \eps) P[N] {\log \log P[N] \over \log P[N]},
$$
provided that $N$ is sufficiently large in terms of $b,k$ and $\eps$. %%y  
This implies Theorem 1.1.    \cqfd 

\medskip

\noindent {\it Proof of Theorem 1.2.}

We assume that $q_1,\ldots,q_s$ are the prime divisors of $N$, thus in
particular we have $s=\omega (N)$. %%y 
If (3.2) is not satisfied, then $\log \log N \le c_7 k$. Otherwise, 
by taking the logarithms of both sides of (3.3) and using (3.1), we get
$$
{\log \log N \over k} \le c_8 +  C_7 \omega(N) +   \log k 
+ \omega (N) \log \log P[N]  + \log \log (k \log P[N]).       
$$
This establishes Theorem 1.2. 

\medskip

\noindent {\it Proof of Theorem 1.3.}  

We argue as in the proof of Theorem 1.1. 
Let $\eps$ be a positive real number with $\eps < \delta_0/3$. Suppose that 
$j$ is sufficiently large and
set $N:=u_j^{(f)}$. It follows from (1.1) that (3.2) holds if $N$ is large enough. 
Then, (3.3) holds with an integer $k$ at most equal to $f(N)$.  
By using (3.10), (3.11) and the Prime Number Theorem,  we get 
$$
\log \log N \leq (1+\eps)
f(N)\Bigl( \log f(N) + P[N] {\log \log P[N] \over \log P[N]} \Bigr).
$$
It then follows from the definition of the positive real number $\delta_0$ that 
$$
\log\log N <(1-\delta_0+ 3 \eps)\log\log N +(1+\eps) f(N) P[N] {\log \log P[N] \over \log P[N]},
$$
hence
$$
(\delta_0 - 3 \eps)  \Psi_f (N)< (1+\eps) P[N] {\log \log P[N] \over \log P[N]}.
$$
Therefore, we have proved (1.2).     
\cqfd

\vskip 5mm

\goodbreak

\centerline{\bf 4. Additional remarks}

\vskip 5mm

\medskip

\noi {\it Remark 4.1.} 
Arguing as Stewart did in \cite{Ste08c}, 
we can derive a lower bound for $Q[u_j^{(k)}]$, where $Q[n]$ denotes the
greatest square-free divisor of a positive integer $n$, similar to the 
lower bound for $Q[u_j^{(3)}]$ given in Theorem 4.1 of \cite{Bu17}.

\medskip

\noi {\it Remark 4.2.} 
Let $a, b$ be integers 
such that $a > b > 1$ and $\gcd(a, b) \ge 2$. Perfect powers in the double 
sequence $(a^m + b^n + 1)_{m, n \ge 1}$ have been considered 
in \cite{BeBuMi12,BeBuMi13,Bu17,CoZa05}.  
The method of the proof of Theorem 1.2 allows us to establish the 
following extension of Theorem 4.3 of \cite{Bu17}. 

\proclaim Theorem 4.1.
Let $k \ge 2$ and
$a_1,  \ldots  , a_k$ be positive integers with $\gcd (a_1, \ldots , a_k) \ge 2$. 
Let $\bfv = (v_j)_{j \ge 1}$ denote the increasing sequence composed of all the 
integers of the form $a_1^{n_1} + \cdots + a_k^{n_k} + 1$, with $n_1, \ldots , n_k \ge 1$. 
Then, for every positive $\eps$, we have
$$
P [v_j] > \Bigl( {1 \over k-1}  - \eps \Bigr) \, 
\log \log v_j \, {\log \log \log v_j \over \log \log \log \log v_j},
$$
when $j$ exceeds some effectively computable constant 
depending only on $a_1,  \ldots  , a_k$, and~$\eps$. 

\medskip

\noi {\it Remark 4.3.} 
Let $n$ be a positive integer. 
Let $S = \{q_1, \ldots , q_s\}$ be a finite, non-empty set of distinct prime numbers.
Write $n = q_1^{r_1} \cdots q_s^{r_s} M$, where  
$r_1, \ldots , r_s$ are non-negative integers and $M$ is an integer 
relatively prime to $q_1 \cdots q_s$. We define the $S$-part $[n]_S$ of $n$ by 
$$
[n]_S := q_1^{r_1} \cdots q_s^{r_s}. 
$$
Theorem 1.1 of \cite{Bu17}
asserts that, for every $k \ge 2$ and every positive real number $\eps$, 
we have 
$$
[u_j^{(k)}]_S <  (u_j^{(k)})^{\eps}, 
$$
for every sufficiently large integer $j$. This implies that
(and is a much stronger statement than)  
the greatest prime factor of $u_j^{(k)}$ tends to infinity as $j$ tends to infinity.
The proof uses the Schmidt Subspace Theorem and it is here essential that $k$ is fixed. 
Moreover, this is an ineffective result. 

The main goal of \cite{Bu17} was to establish an effective improvement of 
the trivial estimate $[u_j^{(3)}]_S \le u_j^{(3)}$ 
of the form $[u_j^{(3)}]_S \le (u_j^{(3)})^{1 - \delta}$, for a small 
positive real number $\delta$ and for $j$ sufficiently large. 
A key tool was a stronger version of Theorem 2.1 in the special case where $|b_n|$ 
is small. 
Unfortunately, for $k > 3$, the method of the proof of 
Theorem 1.2 does not seem to combine well with this stronger version of Theorem 2.1
to get an analogous result. %%y 
We are only able to establish that, for any fixed integer 
$k \ge 4$ and any given positive real number $\eps$, the upper bound
$$
[u_j^{(k)}]_S \le u_j^{(k)} \, \exp \bigl( - (\log u_j^{(k)})^{(1-\eps)/(k-2)} \bigr)
$$
holds for every sufficiently large integer $j$.

\medskip

\noi {\it Remark 4.4.} 
Instead of considering the number of nonzero digits in the representation of an integer 
in an integer base, we can focus on the number of blocks composed of the same digit 
in this representation, a quantity introduced by Blecksmith,
Filaseta, and Nicol \cite{BlFiNi93}; see also \cite{BaTiTi99,BuCiMi13}. A straightforward 
adaptation of our proofs shows that analogous versions of Theorems 1.1 to 1.3 hold 
with `number of nonzero digits' replaced by `number of blocks'. We omit the details.

\medskip

\noi {\it Remark 4.5.} 
In the opposite direction of our results, it does not seem to be easy to confirm the 
existence of arbitrarily large 
integers with few digits in their representation in some integer base and only small 
prime divisors. A construction given in Theorem 6 of \cite{Bu98} and based 
on cyclotomic polynomials shows that there exist an absolute, positive real number $c$ 
and arbitrarily large integers $N$ of the form $2^n + 1$ such that
$$
P[N] \le N^{c / \log \log \log N}.
$$

\vskip 7mm

\noindent 
{\bf Acknowledgements. } 
The second author was supported by JSPS KAKENHI Grant Number 15K17505. 
%The author is very grateful to the referee for a detailed report. 

\vskip 18mm

\centerline{\bf References}

\vskip 7mm

\beginthebibliography{999}

\medskip

\bibitem{BaTiTi99}
G. Barat, R. F. Tichy, and R. Tijdeman, 
{\it Digital blocks in linear numeration systems}.  
In: Number theory in progress, Vol. 2 
(Zakopane-Ko\'scielisko, 1997),  607--631, de Gruyter, Berlin, 1999.

\bibitem{BeBuMi12}
M. A. Bennett, Y. Bugeaud, and M. Mignotte,
{\it Perfect powers with few binary digits
and related Diophantine problems, II},
Math. Proc. Cambridge Philos. Soc. 153 (2012), 525--540.

\bibitem{BeBuMi13}
M. A. Bennett, Y. Bugeaud, and M. Mignotte,
{\it Perfect powers with few binary digits
and related Diophantine problems},
 Ann. Sc. Norm. Super. Pisa Cl. Sci. 12 (2013), 525--540.

\bibitem{BlFiNi93}
R. Blecksmith, M. Filaseta, and C. Nicol, 
{\it A result on the digits of $a\sp n$},
Acta Arith.  64  (1993),  331--339.

\bibitem{Bu98}
Y. Bugeaud, 
{\it Lower bounds for the greatest prime factor of $a x^m+ b y^n$}, 
Acta Math. Inform. Univ. Ostraviensis 6 (1998), 53--57. 

\bibitem{Bu17}
Y. Bugeaud, 
{\it On the digital representation of integers with bounded prime factors}, 
Osaka J. Math.  
To appear. 

\bibitem{BuCiMi13}
Y. Bugeaud, M. Cipu, and M. Mignotte,
{\it On the representation of Fibonacci and Lucas numbers in an integer base}, 
Ann. Math. Qu\'e. 37 (2013), 31--43.

\bibitem{CoZa05}
P. Corvaja and U. Zannier,
{\it $S$-unit points on analytic hypersurfaces}, 
Ann. Sci. \'Ecole Norm. Sup. 38 (2005), 76--92.

\bibitem{Matv00} 
E.\ M.\ Matveev,
{\it An explicit lower bound for a homogeneous rational linear form
in logarithms of algebraic numbers.\ II},
Izv.\ Ross.\ Acad.\ Nauk Ser.\ Mat.\  {64}  (2000),  125--180 (in Russian); 
English translation in Izv.\ Math.\  {64} (2000),  1217--1269.

\bibitem{Ste80}
C. L. Stewart,
{\it On the representation of an integer in two different bases},
J. reine angew. Math.  319  (1980), 63--72.

\bibitem{Ste08c}
C. L. Stewart,
{\it On the greatest square-free factor of terms of a linear recurrence sequence}. 
In: Diophantine equations, 257--264, 
Tata Inst. Fund. Res. Stud. Math., 20, Tata Inst. Fund. Res., Mumbai, 2008.

\bibitem{Yu07}
K. Yu, 
{\it $p$-adic logarithmic forms and group varieties. III}, 
Forum Math. 19 (2007), \ \ 187--280.

\vskip 1cm

\noi Yann Bugeaud   \hfill{ }

\noi Institut de Recherche Math\'ematique Avanc\'ee, U.M.R. 7501

\noi Universit\'e de Strasbourg et C.N.R.S.

\noi 7, rue Ren\'e Descartes

\noi 67084 Strasbourg, FRANCE

\medskip

\noi e-mail : {\tt bugeaud@math.unistra.fr}

\vskip 1cm

\noi Hajime Kaneko  \hfill{}

\noi Institute of Mathematics, University of Tsukuba, 1-1-1

\noi Tennodai, Tsukuba, Ibaraki, 350-0006, JAPAN

\noi Center for Integrated Research in Fundamental Science and Technology (CiRfSE) 

\noi University of Tsukuba, 

\noi Tsukuba, Ibaraki, 305-8571, JAPAN

\medskip

\noi e-mail: {\tt kanekoha@math.tsukuba.ac.jp}

\bye